\documentclass[reqno,11pt]{amsart}
\usepackage{amsmath,amssymb,tabularx,setspace,color,multirow,graphics,graphicx}
\usepackage[margin=3.5cm]{geometry}
\newtheorem{theorem}{Theorem}
\usepackage{adjustbox}

\newtheorem{definition}{Definition}
\newtheorem{remark}{Remark}
\newtheorem{corollary}{Corollary}
\makeatletter
\renewcommand\section{\@startsection {section}{1}{\z@}%
                                   {-3.5ex \@plus -1ex \@minus -.2ex}%
                                   {2.3ex \@plus.2ex}%
                                   {\normalfont\large\bfseries}}
\makeatother

\newtheorem{Corollary}{Corollary}

\newtheorem{Theorem}{Theorem}

\newtheorem{Example}{Example}

\begin{document}
\doublespace
\title[]{Dynamic Cumulative Residual Entropy Generating Function and its properties}
\author[]%
{S\lowercase{mitha} S.$^{\lowercase{a}}$, S\lowercase{udheesh} K. K\lowercase{attumannil}$^{\lowercase{b}}$  \lowercase{and}   S\lowercase{reedevi} E. P.$^{\lowercase{c}}$  \\
 $^{\lowercase{a}}$K E C\lowercase{ollege} M\lowercase{annanam,} K\lowercase{erala}, I\lowercase{ndia},\\
 $^{\lowercase{b}}$I\lowercase{ndian} S\lowercase{tatistical} I\lowercase{nstitute},
  C\lowercase{hennai}, I\lowercase{ndia,}\\
$^{\lowercase{c}}$M\lowercase{aharaja's} C\lowercase{ollege} E\lowercase{ranakulam}, K\lowercase{erala}, I\lowercase{ndia.}}
\thanks{{$^{\dag}$}{Corresponding E-mail: \tt skkattu@isichennai.res.in }}
\maketitle
\vspace{-0.2in}
\begin{abstract}
In this work, we study the properties of  cumulative residual entropy generating function.
We then introduce dynamic cumulative residual entropy generating function (DCREGF).  It is shown that the
DCREGF determines the distribution uniquely. We study some
characterization results using the relationship between DCREGF and hazard rate
and mean residual life function. A new class of life distribution based on decreasing DCREGF is introduced. Finally we develop a test for decreasing DCREGF and study its performance.
\end{abstract}
\vspace{-0.1in}
\section{Introduction}
Entropy is an important concept in the field of information theory and Shannon (1948) was the first who formally introduced it.  To measure
the uncertainty contained in a random variable $X$,  entropy is defined as
\begin{equation}\label{eq1}
  H(X)=-\int_{0}^{\infty}\log f (x)f(x)dx=E(-\log f(X)),
\end{equation} where `log' denotes the natural logarithm. Several measures of entropy have been introduced in the  literature each one  suitable for some specific situations. The widely used  measures of entropy are cumulative residual entropy (CRE) is given by (Rao et al., 2004)
\begin{equation}\label{rent}
     \mathcal{E}(X)=-\int_{0}^{\infty}\bar{F}(x)\log \bar{F}(x)dx,
\end{equation}
where  $\bar{F}(x)=1-F(x)$ is the survival function of $X$.
Di Crescenzo and Longobardi (2009) introduced  cumulative entropy (CE) for estimating the uncertainty in the past  lifetime of a system as
\begin{equation}\label{cent}
     \mathcal{CE}(X)=-\int_{0}^{\infty}F(x)\log F(x)dx.
\end{equation}
The weighted versions of  $\mathcal{E}(X)$ and $\mathcal{CE}(X)$ have been studied in the literature as well.  These are given by (Mirali et al., 2016)
\begin{equation}\label{wcent}
     \mathcal{E}^{w}(X)=-\int_{0}^{\infty}x\bar{F}(x)\log \bar{F}(x) dx.
\end{equation} and (Mirali and Baratpour, 2017)
\begin{equation}\label{wrent}
     \mathcal{CE}^{w}(X)=-\int_{0}^{\infty}xF(x)\log F(x)dx.
\end{equation} For recent development in this area we refer to Kharazmi  and Balakrishnan  (2020, 2021a, 2021b), Hashempour et al. (2022), Kazemi et al. (2022) and Sudheesh et al. (2022).   Among these,  Sudheesh et al. (2022) defined a  generalized cumulative residual entropy and study its properties  They show that cumulative residual entropy,  weighted cumulative residual entropy are special cases of the proposed measure.

The moment generating function (m.g.f) of a probability distribution is a convenient tool for evaluating mean, variance and other moments of a probability distribution. The successive derivative of the m.g.f at a point zero gives the successive moments of the probability distribution provided  these moments exists. In information theory, generating functions has been defined for probability density function to determine
information quantities such as Shannon information, extropy and Kullback-Leibler divergence.
 Golomb (1966) introduced entropy generating of a probability distribution and is given by
$$
B(s)=\int f^{s}(x) d x ; s>1 .
$$
It may be noted that the first derivative of this function at $s=1$, gives the negative of  Shannon's entropy given in (\ref{eq1}).
%

The study of time to event of a subject is of interest in reliability, survival analysis and many other fields. Considering the importance of $B(s)$, in such practical situations we had developed residual entropy generating function using the residual r.v, $X_{t}=X-t \mid X>t$, is given by
$$
B_{t}(s)=\int_{t}^{\infty}\left(\frac{f(x)}{\bar{F}(t)}\right)^{s} d x ; s>1 .
$$
Also $-\frac{d}{d t} B_{t}(s) \mid s=1$ gives residual entropy introduced by Ebrahimi (1996).

Rao et al. (2004) introduced a new measure of uncertainty called cumulative residual entropy (CRE) in a distribution function $F$ and obtained some properties, which is given by
\begin{equation}\label{cre}
  \xi(X)=-\int_{0}^{\infty} \bar{F}(x) \log \bar{F}(x) d x.
\end{equation}Asadi and Zohrevand (2007) have modified the definition of CRE in order to accommodate the current age of the system is given by
\begin{equation}\label{dcre}
  \xi(X ; t)=-\int_{t}^{\infty} \frac{\bar{F}(x)}{\bar{F}(t)} \log \left(\frac{\bar{F}(x)}{\bar{F}(t)}\right) d x .
\end{equation}The measure $\xi(X ; t)$ is known as dynamic cumulative residual entropy function. Obviously $\xi(X ; 0)=\xi(X)$. For more properties and applications of (\ref{cre}) and (\ref{dcre}) one may refer to Rao et al.  (2004). Asadi and Zohrevand (2007), Navarro et al. (2010) and the references therein.

The rest of the paper is organized as follows. In Section 2 study some properties of cumulative residual entropy generating function (CREGF). We also discuss the nonparametric estimation of CREGF.    In section 3, we introduce dynamic cumulative residual entropy generating function. We prove that the proposed measure determines the distribution uniquely. We also study some characterization results using the relationship of DCREGF with hazard rate and mean residual life function. In Section $4$, we develop a test for testing decreasing DCREGF. In Section 5, we conduct a Monte Carlo simulation   study to assess the finite sample performance the proposed test. Some concluding remarks with some open problems are given in Section 6.

\section{Cumulative residual entropy generating function}
Let $X$ be non-negative continuous random variable having distribution function $F(x0$. Let $\bar{F}(x)$ is the survival function of $X$.
The survival function is more useful in lifetime studies. This motivate us to study the properties of CREGF.

\begin{definition}
  Let $X$ be a non-negative random variable having an absolutely continuous   survival function $\bar{F}(x)$, CREGF denoted by $C_{s}(X)$ is defined as
\begin{equation}\label{cregf}
  C_{s}(X)=\int_{0}^{\infty}(\bar{F}(x))^{s} d x ,\,\,s>0.
\end{equation}
\end{definition}
We noted that Kharazmi and Balakrishnan  (2021) discussed the measure in (\ref{cregf}) and obtained some relationship with Gini mean difference. They also discussed relative cumulative residual information generating measure to study the closeness between two survival function. We further  study some properties of  $C_{s}(X)$.   We also discuss the nonparametric estimation of $C_{s}(X)$.

It may be noted that negative of the first derivative of $C_{s}(X)$ at $s=1$  gives the cumulative residual entropy function given in (\ref{cre}). Also note that (prime denote the derivative)
$$
\left.C_{s}^{\prime}(X)\right|_{s=2}=\int_{0}^{\infty}(\bar{F}(x))^{2} \log \bar{F}(x) d x .
$$
The Table 1 gives the expression of $C_{s}(X)$ for some well-known distributions.
\begin{table}[h]
  \centering
   \caption{Expression of $C_{s}(X)$ for some well-known distributions}
\begin{tabular}{|c|l|c|c|c|}
\hline
 & Distribution & $\bar{F}(x)$ & $C_s(X)$ & $\mathcal{CE}(X)$ \\
\hline
$(i)$ & $U(0, a)$ & $1-\frac{x}{a} ,\,>0$ & $\frac{a}{s+1}$ & $\frac{a}{4}$ \\
\hline
$(i i)$ & $\operatorname{GPD}$ & $\left(1+\frac{a x}{b}\right)^{-\left(1+\frac{1}{a}\right)} ,\,x>0, a>-1, b>0$ & $\frac{b}{(a+1) s-a}$ & $b(a+1)$ \\
\hline
$(i i i)$ & Pareto $(k ; \alpha)$ & $\left(\frac{k}{x}\right)^{\alpha} ,\, x \geq k$ & $\frac{k}{\alpha s-1}$ & $\frac{k \alpha}{(\alpha-1)^{2}}$ \\
\hline
$(i v)$ & $\exp (\lambda)$ & $\exp \{-\lambda x\} ,\, x \geq 0$ & $\frac{1}{\lambda s}$ & $\frac{1}{\lambda}$ \\
\hline
$(v)$ & Pareto $I I$ & $\left(1+\frac{x}{a}\right)^{-b} ,\, x \geq 0$ & $\frac{a}{b s-1}$ & $\frac{a b}{(b-1)^{2}}$ \\
\hline
\end{tabular}
\end{table}
Next, we prove some properties of $C_s(X)$.
In the following property, we show that $C_{s}(X)$ is a shift independent measure.
\begin{theorem}
 Let $X$ be continuous nonnegative random variable and  $Y=a X+b$, with $a>0$ and $b \geq 0$, then we have
$$
C_{s}(Y)=a \cdot C_{s}(X).
$$
\end{theorem}
\noindent{\bf Proof:}  The result follows by noting that $\bar{F}_{a X+b}(x)=\bar{F}_{X}\left(\frac{x-b}{a}\right)$ for all $x>b$.

\begin{theorem}
  Suppose $X$ and $Y$ are two random variable admitting the proportional hazard model given by
  $$
\bar{F}_{\theta}^{*}(x)=(\bar{F}(x))^{\theta},\, \theta>0, x>0.
$$
Then the following statements hold;

(a) $C_{s}\left(Y_{\theta}\right)=C_{s \theta}(X)$

(b) $C_{s}(Y)=\theta \cdot C_{s}(\theta X)$.
\end{theorem}

\begin{corollary}
  Let $X$ be  non-negative random variable having  absolutely continuous distribution function $F$. Let $X_{1: n}$ be the first order statistic based on a  random sample  $X_{1}, X_{2}, \ldots, X_{n}$ from  $F$. We have $\bar{F}_{X_{1: n}}(x)=(\bar{F}(x))^{n}$ and so $C_{s}\left(X_{1: n}\right)=C_{n s}(X)$.
\end{corollary}

\begin{Example}
  Suppose $X \sim \exp (\lambda)$, then $C_{s}\left(Y_{\theta}\right)=\frac{1}{\lambda \theta s}, C_{s \theta}(X)=\frac{1}{\lambda \theta s}$ and $C_{s}(X)=\frac{1}{\lambda s}$. So we have $C_{s}\left(Y_{\theta}\right)=C_{s \theta}(X)$ and $C_{s}(X)=\theta \cdot C_{s \theta}(X)$.
\end{Example}
In the next theorem, we give bound of CREGF based on the mean of $X$.
\begin{theorem}
  Let $X$ be a non-negative continuous random variable  with finite  mean $\mu$ and CREGF $C_{s}(X)$ then $C_{s}(X)<\mu$.
\end{theorem}
\noindent{\bf Proof:} We have,
$$
(\bar{F}(x))^{s}<\bar{F}(x),\, s>1 .
$$
Integrating both sides of the above equation with respect to  $x$,
$$
\int_{0}^{\infty}(\bar{F}(x))^{s} d x<\int_{0}^{\infty} \bar{F}(x) d x .
$$
That is,
$$
C_{s}(X)<\mu .
$$

Next, we find an estimator for $C_{s}(X)$. Here we assume that $s$ is positive integer. Let $X_{1: n}$ be the first order statistic based on a  random sample  $X_{1}, \ldots, X_{n}$ from  $F$.  Then we have $\bar{F}_{X_{1: n}}(x)=(\bar{F}(x))^{n}$. For a non-negative random variable $X$, we have $\mu=E(X)=\int_{0}^{\infty}\bar{F}(x)dx$. Hence, we can write
\begin{equation}\label{cregfnew}
  C_{s}(X)=\int_{0}^{\infty}(\bar{F}(x))^{s} d x=E(X_{1: s}).
\end{equation}
Hence a U-statistic based  estimator of $C_{s}(X)$ is given by
\begin{equation}\label{est}
 \widehat{C}_{s}(X) =\frac{1}{C_{s,n}} \sum\limits_{C_{s,n}}\min{(X_{i_1},X_{i_2},\cdots,X_{i_{s}})},
\end{equation}where the summations is over the set  $C_{s,n}$ of all combinations of $s$ distinct elements $\lbrace i_1,i_2,\cdots,i_{s}\rbrace$ chosen from $\lbrace 1,2,\cdots,n\rbrace$.
Now, we express $\widehat{\Delta} $ in a simple form.
Let $X_{i: n}$ be the $i$-th order statistics based on $n$ random sample $X_1,\ldots,X_n$ from $F$.
In terms of the order statistics we have the following equivalent expressions
\begin{equation*}
\sum\limits_{i = 1}^n \sum\limits_{j = 1,j<i}^n\min\{X_1,X_2\} =\sum\limits_{i = 1}^n {(n - i)X_{1: n} }
\end{equation*}
 and
\begin{eqnarray*}
\sum\limits_{i = 1}^n \sum\limits_{j= 1,j<i}^n\sum\limits_{k = 1,k<j}^n\min\{X_1,X_2,X_3\} &=&\sum\limits_{i = 1}^n {\frac{{(n - i - 1)(n - i)X_{1: n} }}{2}}  \\&=& \sum\limits_{i = 1}^n {\binom{n - i}{2} } X_{1: n}.
\end{eqnarray*}
Therefore,  the  estimator given in (\ref{est}) can be expressed  as
\begin{equation}\label{uses}
 \widehat{C}_{s}(X) =\frac{1}{C_{s,n}}  \sum\limits_{i = 1}^n {\binom{n - i}{s-1} X_{1: n} }.
\end{equation}
Next we study the asymptotic properties of  $ \widehat{C}_{s}(X)$. Clearly $ \widehat{C}_{s}(X)$ is an unbiased and consistent estimator of  ${C}_{s}(X)$ (Lehmann, 1951). In the following result we establish the asymptotic distribution of $ \widehat{C}_{s}(X)$.
\begin{theorem}
  As $n\rightarrow \infty$,  $\sqrt{n}(\widehat{C}_{s}(X)-{C}_{s}(X))$ converges in distribution to a normal random variable with mean zero and variance $s^2\sigma^2$, where $\sigma^2$ is given by
  \begin{small}
 \begin{equation}\label{estvar}
\sigma^{2}=Var\Big(X{{\bar F}^{s- 1}}(X) + (s-1)\int_0^X {{{y \bar F}^{s  - 2}}(y)} d{F}(y)\Big).
\end{equation}
\end{small}
\end{theorem}
\noindent {\bf Proof:} Using the central limit theorem for U-statistics  we have the asymptotic  normality of $ \widehat{C}_{s}(X)$. The asymptotic variance is $s^2\sigma_1^2$ where $\sigma_1^2$ is given by (Lee, 2019)
\begin{equation}\label{estvar1}
\sigma_1^2= Var\left[E\left(\min(X_1,\ldots,X_{s})|X_{1}\right)\right].
\end{equation}Denote $Z=\min(X_2,X_3,...,X_s)$, then the survival function of $Z$ is given by $\bar{F}^{s-1}(x)$.  Consider
\begin{eqnarray*}\label{eq4}
E\left[ {\min \left( {x,{X_2},{X_3},...,{X_{s} }} \right)} \right]& =& E\left[ {xI(Z > x)} \right] + E\left[ {ZI(Z \le x)} \right]\nonumber\\&=&x{{\bar F}^{s- 1}}(x) + (s-1)\int_0^x {{{y \bar F}^{s  - 2}}(y)} d{F}(y).
\end{eqnarray*}
Therefore, from (\ref{estvar1}) we obtain the variance expression specified in  equation (\ref{estvar}). THis completes the proof of the theorem.
\\ The finite sample performance of the estimator given in (11) is evaluated through Monte Carlo simulation and the results of the same is reported in Section 5.
\vspace{-0.2in}
\section{Dynamic cumulative residual entropy generating function }\vspace{-0.1in}
In this section, we define dynamic cumulative residual entropy generating function (DCREGF). We  establish some characterization results for some  well-known distributions in terms of DCREGF.
 Let $X$ be the lifetime of a component or system under the condition that the system has survived up to an age $t$. In such cases, we are interested in studying  dynamic or time dependent random variable $X_{t}=(X-t \mid X>t)$ having survival function
$$
\bar{F}_{t}(x)= \begin{cases}\frac{\bar{F}(x)}{\bar{F}(t)}, & \text { where } x>t \\ 1, & \text { otherwise. }\end{cases}
$$
Using the survival function of  $X_{t}$, we define  DCREGF of $X$.

\begin{definition}Let $X$ be a continuous random variable with distribution function $F(x)$.  Then DCREGF of $X$ is defined as
   \begin{equation}\label{dcregf}
  C_{s}(X ; t)=\int_{t}^{\infty}\left(\frac{\bar{F}(x)}{\bar{F}(t)}\right)^{s} d x ,\,\, s \geq 1.
\end{equation}
\end{definition}

From  Definition 2, we readily observe the following properties of $C_s(X)$.
\begin{enumerate}
  \item For $t=0, C_{s}(X ; t)=C_{s}(X)$;
  \item When $s=1$,  $C_{1}(X ; t)=m(t)$;
  \item Differentiating  $C_{s}(X ; t)$ with respect to $s$ and taking $s=1$, we obtain the negative of the dynamic cumulative residual entropy given in (\ref{dcre}).
\end{enumerate}
\begin{theorem}
  For a non negative random variable $X$, let $Y=a X+b$, where $a>0$ and  $b \geq 0$. Then we have  $$C_{s}(Y ; t)=a \cdot C_{s}\left(X ; \frac{t-b}{a}\right),\,\,\, t \geq b.$$
\end{theorem} \noindent{\bf Proof:} Consider
$$
\begin{aligned}
C_{s}(Y ; t) &=\int_{t}^{\infty}\left(\frac{\bar{F}_{X}\left(\frac{x-b}{a}\right)}{\bar{F}(t)}\right)^{s} d x \\
&=a \int_{\frac{t-b}{a}}^{\infty}\left(\frac{\bar{F}_{x}(u)}{\bar{F}(t)}\right)^{s} d u \\
&=a \cdot C_{s}\left(X ; \frac{t-b}{a}\right), t \geq b.
\end{aligned}
$$
From the above theorem, the following property can be easily obtained for  particular choice of $C_{s}(Y ; t)$.
\begin{corollary}We have

  (a) $C_{s}(a X ; t)=a \cdot C_{s}\left(X ; \frac{t}{a}\right)$.

(b) $C_{s}(X+b ; t)=C_{s}(X ; t-b)$.
\end{corollary}
\begin{theorem}Let $X$ be a continuous random variable with distribution function $F(x)$. Also defined the hazard rate of $X$ as $h(t)=f(t)/\bar{F}(t)$.
 The relationship between hazard rate and DCREGF is given by
\begin{equation}\label{rhd}
  h(t) =\frac{1+C_{s}^{\prime}(X ; t)}{s \cdot C_{s}(X ; t)}.
\end{equation}
\end{theorem}
\noindent{\bf Proof:}  From the definition in (\ref{dcregf}), we have
$$
(\bar{F}(t))^{s} C_{s}(X ; t)=\int_{t}^{\infty}(\bar{F}(x))^{s} d x .
$$
Differentiating both side of the above equation  with respect to  $t$, we obtain
$$
\begin{aligned}
(\bar{F}(t))^{s} C_{s}^{\prime}(X ; t)-C_{s}(X ; t) \cdot s \cdot(\bar{F}(t))^{s-1} f(t) &=-(\bar{F}(t))^{s},
\end{aligned}
$$where prime denotes the derivative. In terms of hazard rate, the above eqaution can be written as
$$C_{s}^{\prime}(X ; t)-s \cdot h(t) \cdot C_{s}(X ; t)+1 =0.$$
Hence, we obtain
$$
\begin{aligned}
h(t) &=\frac{1+C_{s}^{\prime}(X ; t)}{s \cdot C_{s}(X ; t)}.
\end{aligned}
$$
In the following theorem, we show that the dynamic cumulative residual entropy generating function determines the distribution of $X$ uniquely.
\begin{theorem}\vspace{-0.2in}
 Let $X$ be a non-negative random variable  with density function $f(x)$, the survival function $\bar{F}(x)$ and the hazard rate $h(x)$. Then $C_{s}(F ; t)$ uniquely determines the survival function $\bar{F}(t)$.
\end{theorem}\vspace{-0.2in}
\noindent{\bf Proof:} From the relationship between the hazard rate and the DCREGF given (\ref{rhd}), we have
$$-\frac{d}{d t} \log \bar{F}(t) =\frac{1+C_{s}^{\prime}(X ; t)}{s \cdot C_{s}(X ; t)}. $$
 Integrating over the interval $(0,t)$, we obtain
 $$ -\log \bar{F}(t) =\int_{0}^{t}\frac{1+C_{s}^{\prime}(X ; t)}{s \cdot C_{s}(X ; t)}dt.$$
That is
$$
\bar{F}(t)=\exp \left[-\int_{0}^{t}\left(\frac{1+C_{s}^{\prime}(X ; t)}{s . C_{s}(X ; t)}\right) d t\right] \text {. }
$$
This shows that the knowledge of $C_{s}(X ; t)$ enables us to determine the distribution.

Now, suppose that $F(x)$ and $G(x)$ are two distribution functions such that
\begin{equation}\label{rbdyn}
  C_{s}(F ; t)=C_{s}(G ; t).
\end{equation}
That is,
$$
\frac{1}{(\bar{F}(t))^{s}} \int_{t}^{\infty}(\bar{F}(x))^{s} d x=\frac{1}{(\bar{G}(t))^{s}} \int_{t}^{\infty}(\bar{G}(x))^{s} d x.
$$
Differentiating with respect to  $t$, we have
$$
\begin{aligned}
&{[\bar{F}(t)]^{-s} \cdot-(\bar{F}(t))^{s}+\left(\int_{t}^{\infty}(\bar{F}(x))^{s} d x\right) \cdot s(\bar{F}(t))^{-s-1} f(t)} \\
&=[\bar{G}(t)]^{-s} \cdot-(\bar{G}(t))^{s}+\left(\int_{t}^{\infty}(\bar{G}(x))^{s} d x\right) \cdot s(\bar{G}(t))^{-s-1} g(t).
\end{aligned}
$$
In terms of hazard functions  $h_1=f(t)/\bar{F}(t)$ and $h_2=g(t)/\bar{G}(t)$, we have
$$
h_{1}(t) \cdot \frac{1}{(\bar{F}(t))^{s}} \int_{t}^{\infty}(\bar{F}(x))^{s} d x=h_{2}(t) \cdot \frac{1}{(\bar{G}(t))^{s}} \int_{t}^{\infty}(\bar{G}(x))^{s} d x,
$$
which gives
$$
h_{1}(t) C_{s}(F ; t)=h_{2}(t) C_{s}(G ; t).
$$
In view of (\ref{rbdyn}), we obtain
$$
h_{1}(t)=h_{2}(t).
$$
This implies that  $C_{s}(F ; t)$   determines the distribution of $X$ uniquely.\\Next theorem shows that the dynamic cumulative residual entropy generating function is independent of  $t$  if and only $X$ is a exponential random variable.
\begin{theorem}\vspace{-0.2in}
  If $X$ is a non-negative random variable admitting an absolutely continuous  distribution function $F(x)$, then the dynamic cumulative residual entropy generating function is independent of  $t$ if and only if $X$ has exponential distribution.
\end{theorem}\vspace{-0.2in}
\noindent {\bf Proof}. Let $C_{s}(F ; t)=k $, where $k$ is a positive constant. Hence
$$\vspace{-0.1in}
C_{s}^{\prime}(F ; t)=0 .
$$\vspace{-0.2in}
From (\ref{rhd}), we obtain
\begin{equation*}
  sk. h(t)=1.
\end{equation*} Hence, $h(t)=\frac{1}{s k}=\beta$,  a constant.
Constant hazard rate  characterizes the exponential distribution. Therefore, $X$ is distributed as an exponential random variable with parameter  $\beta$.

Conversely, assume that $X \sim \exp (\beta)$ where the survival function is given by  $\bar{F}(x)=\exp(-\beta x)$. Therefore,
$$
\begin{aligned}
C_{s}(F ; t) &=\frac{1}{e^{-\beta t s}} \int_{t}^{\infty} e^{-\beta s x} d x \\
&=\frac{1}{\beta s},
\end{aligned}
$$which is a constant. That is, DCREGF is independent  of  $t$ if and only if $X$ has exponential distribution.

\noindent The next theorem provides a characterization result for the generalized Pareto distribution based on a functional form for the dynamic cumulative residual entropy generating function.
\begin{theorem}\vspace{-0.2in}
  Let $X$ be a non-negative continuous random variable with survival function $\bar{F}(x)$. The dynamic cumulative residual entropy generating function is a linear function of $t$ if and only if $X$ follows generalized Pareto distribution.
\end{theorem}\vspace{-0.2in}
\noindent{Proof:} Assume that, $C_{s}(F ; t)=a+b t ; b \neq 0$. Then $C_{s}^{\prime}(F ; t)=b$. From (\ref{rhd}),
$$
\begin{aligned}
b-s \cdot h(t)(a+b t)+1 &=0.
\end{aligned}
$$Or
$$h(t)(a+b t) =\frac{b+1}{s}.$$
Differentiating with respect to   $t$,  we obtain
$$h(t) b+(a+b t) h^{\prime}(t) =0.$$
That is
$$a+b t =-\frac{b h(t)}{h^{\prime}(t)}.$$Or
$$
\begin{aligned}
-\frac{h^{\prime}(t)}{h(t)} &=\frac{b}{a+b t}=\frac{1}{k+t} ; k=\frac{a}{b}.
\end{aligned}
$$
The above equation can be written as
$$
-\frac{d}{d t} \log h(t)=\frac{1}{k+t} .
$$
Integrating with respect to $t$, we have
$$
\begin{aligned}
&-\log h(t)=\log (k+t)+\log c.
\end{aligned}
$$Or
$$-\log h(t)=\log ((k+t) c).$$
That is\vspace{-0.1in}
$$
\begin{aligned}
\frac{1}{h(t)} &=(k+t) c \\
h(t) &=\frac{1}{c t+d} ,
\end{aligned}
$$where $ d=k c.$ This is the hazard rate of GPD. Since the distribution function uniquely determined by the hazard rate, $X \sim G P D$.

Conversely assume that, $X \sim G P D$, where the survival function of $X$ is given by
$$
\bar{F}(x)=\left(1+\frac{a x}{b}\right)^{-\left(1+\frac{1}{a}\right)}.
$$
Using (14) we get
$$
C_{s}(t)=k(b+a t),\quad k=\frac{1}{(a+1) s-a}.
$$
That is, $C_{s}(t)$ is linear function in  $t$. Hence the theorem.

Now, we provide some characterization results in terms of relationship between the DCREGF and the hazard rate function $h(t)$.
\begin{theorem}\vspace{-0.1in}
  Let $X$ be a non-negative random variable with survival function $\bar{F}(t)$, hazard rate function $h(t)$ and DCREGF  $C_{s}(t)$, then  the following relationship
\begin{equation}\label{rel16}
  C_{s}(t)=k(h(t))^{-1},
\end{equation}
where $k$ is a positive constant, holds if and only if $F(x)$ is the GPD with survival function
$$
\bar{F}(x)=\left(1+\frac{a x}{b}\right)^{-\left(1+\frac{1}{a}\right)},\,\, a>-1, b>0.
$$
\end{theorem}\vspace{-0.2in}
\noindent{\bf Proof}: Under the assumption that (\ref{rel16}) holds, we obtain
$$
C_{s}^{\prime}(t)=-k(h(t))^{-2} h^{\prime}(t).
$$
Using (\ref{rhd}), we have
$$
(s k-1)=-k \frac{h^{\prime}(t)}{(h(t))^{2}}.
$$
That is
$$
\frac{d}{d t}\left(\frac{1}{h(t)}\right)=\left(s-\frac{1}{k}\right) .
$$
Integrating on both sides of the above equation, we have
$$
\frac{1}{h(t)}=\left(\frac{s k-1}{k}\right) t+d_{2} .
$$
That is
$$
\frac{1}{h(t)}=d_{1} t+d_{2} .
$$
Therefore
\begin{equation}\label{eq17}
  h(t)=\frac{1}{d_{1} t+d_{2}},
\end{equation}
where $d_{1}=\frac{s k-1}{k}$ and $d_{2}^{-1}=h(0)$.
Hall and Wellner (1981) has shown that  (\ref{eq17}) is a characteristic property of GPD. Conversely, assume that the random variable $X$ follows GPD. By direct calculation, we obtain
$$
\begin{aligned}
C_{s}(t) &=\frac{b+a t}{(a+1) s-a} \\
&=\left(\frac{b+a t}{a+1}\right)\left(\frac{a+1}{(a+1) s-a}\right) \\
&=k \frac{1}{h(t)}=k(h(t))^{-1}.
\end{aligned}
$$where $ k=\frac{a+1}{(a+1) s-a}$.
Hence we have the proof of the theorem.\\
In the following theorem, we characterize the distribution of $X$ using the relationship between the DCREGF and the mean residual life function (MRL).
\begin{theorem}\vspace{-0.2in}
   If $X$ is a non-negative random variable admitting an absolutely continuous distribution function  $F(x)$. Define MRL of $X$ as $m(t)=E(X-t|X>t)$. Then the relationship, $C_{s}(t)=k . m(t)$  holds for every $t>0$, if and only if $X \sim G P D$.
\end{theorem}\vspace{-0.2in}
\noindent{\bf Proof}:  Assume
\begin{equation}\label{eq18}
C_{s}(t)=k \cdot m(t).
\end{equation}
Differentiating with respect to  $t$, we have
$$
C_{s}^{\prime}(t)=k \cdot m^{\prime}(t) .
$$
From (15), we have
$$
\text { s. } h(t) \cdot C_{s}(t)-1=k \cdot m^{\prime}(t) .
$$
Using (\ref{eq18}), we obtain
$$
k \cdot s \cdot h(t) \cdot m(t)-1=k \cdot m^{\prime}(t)
$$
Now, we have the relationship between the hazard rate and the mean residual life given by
$$
\frac{1+m^{\prime}(t)}{m(t)}=h(t) .
$$
Therefore, we have
$$
\begin{aligned}
k \cdot s\left(1+m^{\prime}(t)\right) &=1+k \cdot m^{\prime}(t) \\
(k s-k) m^{\prime}(t) &=1-k \cdot s \\
m^{\prime}(t) &=\frac{1-k s}{k s-k}=\text { a constant},
\end{aligned}
$$
which implies that $m(t)$ is a linear function in $t$. Linear mean residual life function  is a characteristic property of GPD.

Conversely, assume that $X \sim G P D$. By direct calculation, we obtain
\begin{equation}\label{eq19}
  C_{s}(X ; t)=k(b+a t)=k \cdot m(t),
\end{equation}
where $k=\frac{1}{(a+1) s-a}$. This proves the if part.
\begin{remark}
  Differentiate both sides of (\ref{eq19}) with respect to  $s$ and setting $s=1$ and evaluating  it with negative sign, we obtain Theorem 4.8 of Asadi and Zohervand (2007).
\end{remark}

%

\noindent Next, using DCREGF we introduce two new classes of lifetime distributions.
%
\begin{definition}
   A random variable $X$ is said to have increasing (decreasing) DCREGF, denoted by IDCREGF (DDCREGF) if $C_{s}(X ; t)$ is increasing (decreasing) function in $t, \forall t \geq 0$. 
\end{definition}
\noindent The following theorem gives the bounds for $C_{s}(X ; t)$ in terms of hazard rate function.
\begin{theorem}\vspace{-0.2in}
  The distribution function  $F$ is increasing (decreasing) DCREGF if and only if for all $t>0$
$$
C_{s}(F ; t) \geq(\leq) \frac{1}{s \cdot h(t)},\,\,\,s \geq 1.
$$
\end{theorem}\vspace{-0.2in}
\noindent{\bf Proof:} From the definition of new classes of lifetime distributions, the distribution function  $F$ is said to be increasing (decreasing) if $\left.C_{s} F ; t\right)$ is increasing (decreasing) in $t$. That is \vspace{-0.1in}$$C_{s}^{\prime}(F ; t)\ge(\le) 0.$$ Hence
$$\vspace{-0.1in}
\begin{gathered}
s \cdot h(t) \cdot C_{s}(F ; t)-1\ge(\le) 0.
\end{gathered}\vspace{-0.1in}
$$Or
$$C_{s}(F ; t)\ge (\le) \frac{1}{s \cdot h(t)}.$$
\begin{theorem}\vspace{-0.2in}
  Let $X$ and $Y$ be two non-negative absolutely continuous random variable's with survival functions $\bar{F}(t)$ and $\bar{G}(t)$ and hazard rate functions $h_{1}(t)$ and $h_{2}(t)$, respectively. If $X \geq^{h r} Y$, that is, $\left.h_1(t\right) \leq h_{2}(t)\, \forall t \geq 0$, then $C_{s}(F ; t) \geq$ $C_{s}(G ; t)$.
\end{theorem}\vspace{-0.2in}
\noindent{\bf Proof:} Let  $X \geq^{h r} Y$, then we have
$$
\begin{aligned}
h_{1}(t) \leq h_{2}(t) & \Longrightarrow \quad \bar{F}_{X_{t}}(t) \geq \bar{G}_{Y_{t}}(t) \\
& \Longrightarrow \frac{\bar{F}(x)}{\bar{F}(t)} \geq \frac{\bar{G}(x)}{\bar{G}(t)} \\
& \Longrightarrow \int_{t}^{\infty}\left(\frac{\bar{F}(x)}{\bar{F}(t)}\right)^{s} d x \geq\left(\frac{\bar{G}(x)}{\bar{G}(t)}\right)^{s} \\
& \Longrightarrow C_{s}(F ; t) \geq C_{s}(G ; t).
\end{aligned}
$$

\noindent The following theorem shows that the exponential distribution is the only distribution which is both IDCREGF and DDCREGF.
\noindent The following example gives an application of Theorem 11 in order statistics.
\begin{Example}
  Let $X_{1},  \ldots, X_{n}$ be independent and identical non-negative random variable's with survival function $\bar{F}(x)$. If $X_{i: n}$ denotes the $i$- th order statistic based on a random sample $X_{1}, \ldots, X_{n}$ from $F$, then the following results holds.

(i) $C_{s}\left(X_{i: n} ; t\right) \leq C_{s}\left(X_{i+1: n} ; t\right)$

(ii) $C_{s}\left(X_{i: n} ; t\right) \leq C_{s}\left(X_{1: n-1} ; t\right)$

(iii) $C_{s}\left(X_{n: n} ; t\right) \geq C_{s}\left(X_{n-1: n-1} ; t\right)$.
\end{Example}
\begin{theorem}
 Let $X$ be a non-negative random variable having IDCREGF and DDCREGF, then $X$ has exponential distribution.
\end{theorem}%
\noindent{\bf Proof:} As $X$ possess both IDCREGF and DDCREGF properties, $C_{s}(X ; t)$ is a constant. Hence, from Theorem 8, $X$ has exponential distribution.
\begin{corollary}
Let $X$ be $D D C R E G F$ (IDCREGF), then
$$
\bar{F}(t) \geq(\leq) \exp \left[-\int_{0}^{t} \frac{1}{s . C_{s}(F ; t)} d x\right].
$$
\end{corollary}
\noindent{\bf  Proof:}  Let $X$ be DDCREGF, then
$$
h(t) \leq \frac{1}{s \cdot C_{s}(F ; t)}.
$$
Again, using the relation relationship between $\bar{F}(t)$ and $h(t)$, we have
$$
\bar{F}(t)=\exp \left[-\int_{0}^{t} h(u) d u\right] \geq \exp \left[-\int_{0}^{t} \frac{1}{s \cdot C_{s}(F ; t)} d x\right].
$$
Similarly, we can show that if $X$ has IDCREGF, then
$$
\bar{F}(t) \leq \exp \left\{-\int_{0}^{t} \frac{1}{s \cdot C_{s}(F ; t)} d x\right\}.
$$
\section{Test for decreasing DCREGF}\vspace{-0.5in}
In the previous section, we proved that the DCREGF characterize the distribution of $X$.  As constant DCREGF is a characterization  of  exponential random variable,  we develop a test for testing exponentiality  against the decreasing DCREGF class.

Let $X_1,\ldots X_n$ be a random sample of size $n$ from $F$. We are interested in testing the null hypothesis\vspace{-0.1in}
$$H_0: X \text{ has exponential distribution} $$
 against the alternative\vspace{-0.1in}
 $$H_1: X \text{ has decreasing DCREGF and not exponential.}$$
  For testing the above hypothesis first we define a departure measure which discriminate between null and alternative hypothesis.
  Note that $C_{s}(X ; t)$ is decreasing in $t$ if $C_{s}^{'}(X ; t)\le 0.$ That is,
  $$\frac{sf(t)}{\bar{F}^{s+1}(t)}\int_{t}^{\infty}\bar{F}^s(x)dx-\frac{1}{\bar{F}^{s}(t)}\bar{F}^{s}(t)\le 0$$
  or
  $${\bar{F}^{s+1}(t)}-{sf(t)}\int_{t}^{\infty}\bar{F}^s(x)dx\ge 0.$$
  Hence, we consider  a measure of departure $ \Delta(F)$ given by
\begin{eqnarray}\label{deltam}
\Delta(F)=\int_{0}^{\infty}\left({\bar{F}^{s+1}(t)}-{sf(t)}\int_{t}^{\infty}\bar{F}^s(x)dx\right)dt.
\end{eqnarray}Clearly, $\Delta(F)$ is zero under  $H_0$ and positive  under $H_1$. Accordingly,  $\Delta(F)$  can be considered as a measure of departure  from $H_0$ towards  $H_1$. As the proposed test is based on U-statistics, first we express  $\Delta(F)$ in terms of expectation of the function of random variables. For a non-negative random variable $X$, $E(X)=\int_{0}^{\infty}\bar{F}(x)dx$. Also, observe that $\bar{F}^{n}(x)$ is the survival function of $\min(X_1,\ldots,X_n)$.  Consider
\begin{eqnarray}\label{delta}
\Delta(F)&=&\int_{0}^{\infty}\left({\bar{F}^{s+1}(t)}-{sf(t)}\int_{t}^{\infty}\bar{F}^s(x)dx\right)dt\nonumber\\
&=&E(\min(X_1,\ldots,X_{s+1}))-\int_{0}^{\infty}\int_{t}^{\infty}{sf(t)}\int_{t}^{\infty}\bar{F}^s(x)dxdt.
\end{eqnarray}
Changing the order of integration, we have
\begin{eqnarray}\label{delta1}
\Delta(F)&=&E(\min(X_1,\ldots,X_{s+1}))-s\int_{0}^{\infty}\bar{F}^s(x)\int_{0}^{t}f(t)dtdx\nonumber\\
&=&E(\min(X_1,\ldots,X_{s+1}))-s\int_{0}^{\infty}\bar{F}^s(x)(1-\bar{F}(x))dx\nonumber\\
&=&(s+1)E(\min(X_1,\ldots,X_{s+1}))-sE(\min(X_1,\ldots,X_{s})).
\end{eqnarray}
We find the test statistic using theory of U-statistics.  Consider a symmetric kernel
\begin{eqnarray*}
  h_1(X_1,\ldots,X_{s+1})= (s+1)\min(X_1,\ldots,X_{s+1})-\frac{1}{s+1}\sum_{C_s}s\min(X_{i_{1}},\ldots,X_{i_{s}}),
\end{eqnarray*}where summation is over the set $C_s$ of all  combination of  $s$ integers $i_1<i_2<\ldots<i_s$ chosen from the set $(1,\ldots,s+1)$.
Hence a U-statistic based  test statistic is given by
\begin{equation}\label{test}
 \widehat{\Delta} =\frac{1}{C_{m,n}} \sum\limits_{C_{m,n}}h{(X_{i_1},X_{i_2},\cdots,X_{i_{s+1}})},
\end{equation}where the summations is over the set  $C_{m,n}$ of all combinations of $(s+1)$ distinct elements $\lbrace i_1,i_2,\cdots,i_{s+1}\rbrace$ chosen from $\lbrace 1,2,\cdots,n\rbrace$.
We reject the null hypothesis $H_0$ against the alternative  $H_1$ for large value of $\widehat{\Delta}$.
\begin{remark}
  When $s=1$, the testing problem reduces to testing decreasing mean residual life, a important problem in the lifetime data analysis.
\end{remark}

We find a critical region of the test using the asymptotic distribution of $ \widehat{\Delta}$.  Next we find the asymptotic distribution of $\widehat{\Delta} $.
\begin{theorem}\vspace{-0.3in}
  As $n\rightarrow \infty$,  $\sqrt{n}(\widehat{\Delta}-\Delta(F))$ converges in distribution to normal random variable with mean zero and variance $(s+1)^2\sigma^2$, where $\sigma^2$ is given by
 \begin{eqnarray}\label{vart}
 \sigma^2&=&Var\Big((s+1)X{{\bar F}^{s}}(X) + s(s+1)\int_0^X {{{y \bar F}^{s -1}}(y)} d{F}(y)\nonumber\\
 &&\quad-\frac{s^2}{(s+1)}X{{\bar F}^{s-1}}(X) -\frac{(s-1)s^2}{(s+1)}\int_0^X {{{y \bar F}^{s -2}}(y)} d{F}(y)\Big).
\end{eqnarray}
\end{theorem}
\noindent {\bf Proof:} Using the central limit theorem for U-statistics  we have the asymptotic  normality of $ \widehat{\Delta}^*$. The asymptotic variance is $(s+1)^2\sigma_1^2$ where $\sigma_1^2$ is given by (Lee, 2019)
\begin{equation}\label{var1}
\sigma_1^2= Var\left[E\left(h(X_{1},\ldots,X_{(s+1)})|X_{1}\right)\right].
\end{equation}Denote $Z=\min(X_2,X_3,...,X_s)$, then the distribution of $Z$ is given by $1-\bar{F}^{s-1}(x)$, where $\bar{F}(x)=1-F(x)$.  Consider
\begin{eqnarray*}\label{eq4}
E\left[ {\min \left( {x,{X_2},{X_3},...,{X_{s} }} \right)} \right]& =& E\left[ {xI(Z > x)} \right] + E\left[ {ZI(Z \le x)} \right]\nonumber\\&=&x{{\bar F}^{s- 1}}(x) + (s-1)\int_0^x {{{y \bar F}^{s  - 2}}(y)} d{F}(y).
\end{eqnarray*}\vspace{-0.1in}Similar way, we obtain
\begin{eqnarray*}\label{eq4}
E\left[ {\min \left( {x,{X_2},{X_3},...,{X_{s+1} }} \right)} \right]&=&x{{\bar F}^{s}}(x) + s\int_0^x {{{y \bar F}^{s  - 1}}(y)} d{F}(y).
\end{eqnarray*}\vspace{-0.1in}Hence \vspace{-0.1in}
\begin{eqnarray*}\label{eq4}
 &&\hskip-1.inE\left(h\big( {{X_1},{X_2},...,{X_s }}|X_1=x\big) \right)\\&=&(s+1)x{{\bar F}^{s}}(x) + s(s+1)\int_0^x {{{y \bar F}^{s -1}}(y)} d{F}(y)\\
 &&\quad-\frac{s^2}{(s+1)}x{{\bar F}^{s-1}}(x) -\frac{(s-1)s^2}{(s+1)}\int_0^x {{{y \bar F}^{s -2}}(y)} d{F}(y)+\frac{sk}{(s+1)},
\end{eqnarray*}where $k=E(\min(X_2,\ldots, X_{s+1}))$, a constant.
\noindent Therefore, from (\ref{var1}) we obtain the variance expression specified in the theorem.

Under the null hypothesis $H_0$, $\Delta{(F)}=0$. Hence we have the following corollary.
\begin{corollary}\vspace{-0.3in}
 Under $H_0$, as $n\rightarrow \infty$,  $\sqrt{n}\widehat{\Delta}$ converges in distribution to a Gaussian random variable with mean zero and variance $\sigma_0^2$, where $\sigma_0^2$ is given by
 \begin{equation}\label{nullvarg}
   \sigma_0^2= \frac{s}{(4s^2-1)\lambda^2}.
 \end{equation}
\end{corollary}
\noindent{Proof:} For the exponential distribution the mean is equal to the mean residual life function. Hence
\begin{equation*}
  \int_0^x {{{sy \bar F}^{s -1}}(y)} d{F}(y)=\frac{1}{s\lambda}-\bar{F}^s(x)\left(x+\frac{1}{s\lambda}\right).
\end{equation*}Also
\begin{equation*}
  \int_0^x {{{(s-1)y \bar F}^{s -2}}(y)} d{F}(y)=\frac{1}{(s-1)\lambda}-\bar{F}^{s-1}(x)\left(x+\frac{1}{(s-1)\lambda}\right).
\end{equation*}Therefore, the variance expression in (\ref{vart}) reduces to
\begin{eqnarray*}
  \sigma_{1}^{2} &=& Var\left(\frac{s^2}{(s+1)(s-1)\lambda}\bar{F}^{s-1}(X)-\frac{(s+1)}{s\lambda }\bar{F}^{s}(X)\right) \\
   &=& \frac{1}{(s+1)^2\lambda^2}Var\left( \frac{s^2}{(s-1)}\bar{F}^{s-1}(X)-\frac{(s+1)^2}{s} \bar{F}^{s}(X)\right)
   \\
   &=& \frac{1}{(s+1)^2\lambda^2}\Bigg\{Var\left( \frac{s^2}{(s-1)}\bar{F}^{s-1}(X)\right)+Var\left(\frac{(s+1)^2}{s} \bar{F}^{s}(X)\right)\\&&\quad-\frac{2s(s+1)^2}{(s-1)}Cov(\bar{F}^{s}(X),\bar{F}^{s-1}(X))\Bigg\}.
   \\
   &=& \frac{1}{(s+1)^2\lambda^2}\Bigg\{\frac{s^2}{2s-1}+\frac{(s+1)^2}{2s+1}-(s+1)\Bigg\}= \frac{s}{(s+1)^2(4s^2-1)\lambda^2}.
\end{eqnarray*}  Hence the asymptotic null variance is equal to $s/(4s^2-1)\lambda^2$, which prove the theorem.

In view of the null  variance specified in (\ref{nullvarg}), we consider a scale invariant test given by
$$\widehat{\Delta}^*=\frac{\widehat{\Delta}}{\bar{X}},$$where $\bar{X}$ is the sample mean.  Using Slutsky's theorem we have the following result.
\begin{corollary}
 Under $H_0$, as $n\rightarrow \infty$,  $\sqrt{n}\widehat{\Delta}^*$ converges in distribution to a Gaussian random variable with mean zero and variance $\sigma_0^2$, where $\sigma_0^2$ is given by
 \begin{equation}\label{nullvar}
   \sigma_0^2= \frac{s}{(4s^2-1)}.
 \end{equation}
\end{corollary}

An asymptotic critical region of the test can be obtained  using Corollary 5.  We reject the null hypothesis $H_{0}$ against the alternative hypothesis $H_{1}$ at a significance level $\alpha$, if
\begin{equation*}
 \frac{ \sqrt{n(4s^2-1)} |\widehat{\Delta}^*| }{\sqrt{s}}>Z_{\alpha/2},
  \end{equation*}
where $Z_{\alpha}$ is the upper $\alpha$-percentile point of the standard normal distribution.


\section{Simulation and data analysis}
We conduct Monte Carlo simulation studies using R software to evaluate the finite sample performance of the estimator $\widehat {C}_s(X)$ and the proposed test for DDCREGF. The simulation is repeated ten thousand times.

First, we evaluate the performance of the estimator $\widehat {C}_s(X)$ in terms of bias and MSE. In the simulation study, we generate observations from various lifetime distributions including exponential, gamma, Weibull, lognrmal and Makeham distributions. Various parameters are chosen for generating these observations. The MSE and the absolute bias of the estimators based on samples of sizes $n=10,20,30,40$ and $50$ are calculated. Different choices of $s$ are used in the simulation study. The results of the simulation study are presented in Table 2.  From Table 2, we can observe that both the MSE and the absolute bias are negligible for all distributions and both decrease as $n$ increases.

Next, we conduct an extensive simulation study to assess the performance of the proposed test for decreasing DCREGF. The exponential distribution with different choices of parameters is used to find the empirical type I error of the test.  For finding the empirical power of the test, lifetime distributions including, gamma, Weibull, lognormal, Makeham and linear failure rate, which are members of the DDCREGF class are used.  Random samples of sizes $n=10,20,30,40$ and $50$ are generated from these distributions where the parameters are chosen in such a way that the distribution belongs to decreasing DCREGF family. Different choices of the  $s$ are considered in the study. The results of the simulation study are presented in Tables 3-5. In Table  3, 4 and 5 we reported the results for $s=1$, $s=2$ and $s=3$, respectively.

From Tables 3-5, we can see that empirical type I error of the test approaches chosen significance level for all choices of parameters of the exponential distribution.  For all choices of alternatives, the test yields good power also. When the samples are generated from a gamma distribution with parameters such that, the distribution is approaching exponential, we observe small power. The random samples from Makeham distribution also show low power, compared to the other distributions. The proposed test yield very good power for all the other distributions with various choices of parameters we considered in the simulation.

Next, we use two real life data sets to illustrate the proposed testing procedure. \\
\textbf{Example 1.}  The following data set represents the failure times (in minutes) for a sample of 15 electronic components in an accelerated life test, Lawless (2011).\\
Data: 1.4, 5.1, 6.3, 10.8, 12.1, 18.5, 19.7, 22.2, 23.0, 30.6, 37.3, 46.3, 53.9, 59.8 ,66.2.\\
After applying the above testing procedure, we obtain the test statistic value (for $s=1$) as $0.8409$. Hence,  we accept the null hypothesis that above data set  is exponentially distributed. Our conclusion is same as previous studies on this data.    \\
\textbf{Example 2.}  The second data set is also taken from Lawless (2011).  We considered the data which arose in tests on endurance of deep groove ball bearings. The data are the number of million
revolutions before failure for each of the 23 ball bearings in the life tests and the complete data is given below.\\
Data :17.88, 28.92, 33.00, 41.52, 42.12, 45.60, 48.40, 51.84, 51.96, 54.12, 55.56, 67.80,\\
68.64, 68.64, 68.88, 84.12, 93.12, 98.64, 105.12, 105.84, 127.92, 128.04, 173.40.\\
We obtain the test statistic value (for $s=1$)  as $3.5099$. This concludes that for this data set, we reject the null hypothesis and that it is DCREGF. We use the the function `fitdist' in R-package `fitdistrplus' to fit the ball bearing data. Weibull and lognormal distributions are the two possible models which fit the data and belong to  DCREGF class of distributions.

\section{Conclusions}

In this article, we studied the properties of the cumulative residual entropy generating function (CREGF). We obtain the non-parametric estimator of CREGF and evaluated its performance. We then introduced a dynamic version of the cumulative residual entropy generating function (DCREGF).  We have shown that DCREGF determines the distribution uniquely. Further, we studied the relationships between DCREGF, hazard rate and mean residual life function. This led to the development of a non-parametric test for decreasing DCREGF. We evaluated the finite sample behavior of the proposed test through an extensive Monte Carlo simulation study. The empirical type I of the test is well maintained.   The practical applications of the test are illustrated using real-life examples.  The proposed test can be used for testing decreasing mean residual life function.

Different extensions of measures of entropy and extropy are extensively studied in the literature. We can extend the generating function approach to study these measures.  One can consider developing the empirical likelihood and jackknife empirical likelihood inference for these measures.

\begin{table}[h]
\caption{Bias and MSE of the proposed estimator for various distributions }
\begin{adjustbox}{width=\textwidth}
\begin{tabular}{ccccccccccc}\hline
\cline{1-11}
$n$ &   Bias  &  MSE&   Bias  & MSE &   Bias  & MSE & Bias & MSE & Bias  & MSE \\ \hline
\cline{1-11}
 &  \multicolumn{2}{c}{$exp(1)$}  &  \multicolumn{2}{c}{$gamma(2,2)$}  &  \multicolumn{2}{c}{$Weibull(3,1)$}  &  \multicolumn{2}{c}{$lognormal(0.5,0.5)$}  &  \multicolumn{2}{c}{$Makeham(1,1)$} \\ \hline
\cline{1-11}
\multirow{2}{*}{} & \multicolumn{10}{c}{$s=1$}\\ 
\cline{1-11}

10& 0.0026&	0.1014&	0.0015&	0.0495&	0.0009&	0.0106&	0.0043&	0.0988	&		0.0038&	0.0180 \\
 20& 0.0024&	0.0501&	0.0011&	0.0249&	0.0007&	0.0053&	0.0021&	0.0494	&		0.0018&	0.0088  \\
30 &0.0018&	0.0336&0.0006&	0.0171&0.0004&0.0035&0.0014	&0.0327	&		0.0010&	0.0059  \\
40 & 0.0018	& 0.0257& 0.0002& 0.0127& 0.0002& 0.0026& 0.0013& 0.0246& 		0.0009& 0.0044  \\
50 & 0.0006	&0.0196	&0.0001&0.0098	&0.0000	&0.0021	&0.0010&0.0201	&		0.0008&	0.0035  \\
 \hline
 \multirow{2}{*}{} & \multicolumn{10}{c}{$s=2$}\\ 
\cline{1-11}
10& 0.0017&	0.0351&	0.0015&	0.0255&	0.0007&0.0116&0.0037&0.0490	&		0.0026&0.0124
 \\
20& 0.0016& 0.0167& 0.0009& 0.0127& 0.0005& 0.0056& 0.0014& 0.0240& 	0.0005& 0.0061  \\
30 &0.0008&	0.0115&	0.0005&0.0084&0.0004&0.0037&0.0012&0.0159&			0.0004&0.0041   \\
40 & 0.0005&0.0082	&0.0002	&0.0063	&0.0003	&0.0027&0.0011&0.0122&		0.0002&	0.0030   \\
50 & 0.0003	&0.0068	&0.0001&0.0049&0.0001&0.0022&0.0004	&0.0098	&		0.0001&0.0024  \\
 \hline
 \multirow{2}{*}{} & \multicolumn{10}{c}{$s=3$}\\ 
\cline{1-11}
10& 0.0010& 0.0217& 0.0007&0.0215& 0.0004& 0.0125& 0.0021& 0.0417& 		0.0009&0.0108    \\
20& 0.0007& 0.0106& 0.0007& 0.0101& 0.0004& 0.0061& 0.0021& 0.0204& 	0.0003& 0.0049   \\
30 &0.0005&	0.0068&0.0005&0.0068&0.0003	&0.0041	&0.0008	&0.0137	&		0.0002&0.0032  	 \\
40 & 0.0004&0.0053	&0.0003	&0.0049	&0.0001&0.0031&	0.0002&0.0101&		0.0001&	0.0024  \\
50 & 0.0002	&0.0041	&0.0002	&0.0040&0.0001&0.0024&0.0001&0.0078&		0.0001&0.0020   \\
 \hline
 \multirow{2}{*}{} & \multicolumn{10}{c}{$s=4$}\\ 
\cline{1-11}

10& 0.0015&	0.0170&	0.0011&	0.0189&	0.0011&0.0141&0.0016&0.0403	&		0.0005&	0.0092 \\
20& 0.0005& 0.0077& 0.0006& 0.0092& 0.0011& 0.0067& 0.0013& 0.0191& 	0.0004& 0.0042  \\
30& 0.0005&0.0050&0.0004&0.0058	&0.0007&0.0044&	0.0011&0.0126&		0.0002&0.0028 \\
40 &0.0002&0.0038&0.0003&0.0044	&0.0003	&0.0032	&0.0006&0.0093&		0.0001&	0.0020  \\
50 &0.0001&0.0030&0.0001&0.0034	&0.0001&0.0027&0.0001&0.0075&		0.0001&	0.0016 \\

 \hline
 \multirow{2}{*}{} & \multicolumn{10}{c}{$s=5$}\\ 
\cline{1-11}
10& 0.0010	&0.0136	&0.0009	&0.0177&0.0011&0.0153&	0.0026&	0.0411	&		0.0007&0.0081    \\
20& 0.0004&0.0062&0.0005&0.0082	&0.0008	&0.0073	&0.0022&0.0193&		0.0006&0.0037  \\
30& 0.0004&0.0039&0.0004&0.0054	&0.0004&0.0048&0.0021&0.0126&
0.0004	&0.0025  \\
40 & 0.0004	& 0.0029& 0.0002& 0.0040& 0.0004& 0.0035& 0.0013& 0.0094& 	0.0004& 0.0017  \\
50 & 0.0002& 0.0022	& 0.0001& 0.0032& 0.0001& 0.0029& 0.0003& 0.0074& 	0.0001&0.0014   \\

 \hline
\end{tabular}
\end{adjustbox}
\end{table}

\begin{table}[h]
\caption{Power of the proposed test for various distributions-I }
\begin{adjustbox}{width=\textwidth}
\begin{tabular}{cccccccccccccccccc}\hline
\multirow{2}{*}{} & \multicolumn{12}{c}{$s=1$}\\ 
\cline{1-13}
$n/\alpha$ &  $0.01$  &  $0.05$ &  $0.01$  &  $0.05$ &  $0.01$  &  $0.05$ &  $0.01$  &  $0.05$ &  $0.01$  &  $0.05$ &  $0.01$  &  $0.05$  \\ \hline
 &  \multicolumn{2}{c}{$exp(1)$}  &  \multicolumn{2}{c}{$gamma(2,1)$}  &  \multicolumn{2}{c}{$Weibull(2,1)$}  &  \multicolumn{2}{c}{$lognormal(1,0.5)$}  &  \multicolumn{2}{c}{$Makeham(1,1)$} &  \multicolumn{2}{c}{$LFR(1)$}
 \\ \cline{1-13}

10& 0.0126&0.0605&0.0832&0.2554&0.3379&0.6540&0.4643&0.7713&
0.0556	&0.1811	&0.3385&0.6508    \\
20& 0.0114&	0.0549&	0.2329&	0.4919&0.8017&	0.9531	&0.9021&	0.9772 &0.1346&	0.3391&0.8117&	0.9584   \\
30& 0.0108&0.0544&0.4118&0.6874	&0.9709	&0.9970	&0.9864	&0.9987
&0.2458&0.4965&	0.9721&0.9964  \\
40 & 0.0105&0.0528	&0.5821&0.8276	&0.9966&0.9998&0.9989&0.9999&
0.3678&	0.6382&0.9964&0.9999  \\
50& 0.0099&	0.0510&	0.7232&0.9062&0.9996&1.0000&0.9999&	1.0000
&0.4919&0.7523&0.9995&1.0000  \\
 \hline
 &  \multicolumn{2}{c}{$exp(2)$}  &  \multicolumn{2}{c}{$gamma(2,2)$}  &  \multicolumn{2}{c}{$Weibull(2,2)$}  &  \multicolumn{2}{c}{$lognormal(2,0.5)$}  &  \multicolumn{2}{c}{$Makeham(2,2)$} &  \multicolumn{2}{c}{$LFR(2)$} \\

\cline{1-13}
10&  0.0126&0.0565&	0.0825&0.2537&0.3392&0.6584&0.4709&	0.7768
&0.0540&0.1767&0.3399&0.6653  \\
20& 0.0122&0.0545&0.2289&0.4968	&0.8057&0.9556&0.9782&0.9881
&0.1380	&0.3380	&0.8147&0.9572    \\
30& 0.0121&0.0535&0.4137&0.6929&0.9718&	0.9974&0.9981&0.9998
&0.2405	&0.5004	&0.9715&0.9971 \\
40 & 0.0111&0.0527&0.5746&0.8242&0.9974&1.0000&	1.0000&1.0000
&	0.3723
&0.6477&0.9980&1.0000  \\
50& 0.0098&0.0512&0.7216&0.8995	&0.9999&1.0000&	1.0000&	1.0000
&0.4969&0.7602&0.9999&1.0000 \\
 \hline
 &  \multicolumn{2}{c}{$exp(0.5)$}  &  \multicolumn{2}{c}{$gamma(3,1)$}  &  \multicolumn{2}{c}{$Weibull(3,1)$}  &  \multicolumn{2}{c}{$lognormal(3,0.5)$}  &  \multicolumn{2}{c}{$Makeham(2,1)$} &  \multicolumn{2}{c}{$LFR(0.5)$} \\

\cline{1-13}
10&  0.0129	&0.0608	&0.2513&0.5547&	0.8648&	0.9773&0.4662&0.7695
&0.1195&0.3205&0.3352&0.6546    \\
20& 0.0118&	0.0576&0.6821&0.8987&0.9997	&1.0000	&0.9073	&0.9807
&0.3427&0.6116&0.8097&0.9557    \\
30& 0.0115&0.0533&0.9074&0.9822	&1.0000&1.0000&0.9873&0.9988
&0.5711	&0.8153	&0.9700&0.9967  \\
40 &0.0108&0.0512&0.9801&0.9972	&1.0000	&1.0000	&0.9985&1.0000
&0.7614&0.9227&	0.9978&1.0000   \\
50& 0.0105&0.0495&0.9962&0.9995	&1.0000&1.0000&	0.9999&	1.0000
&0.8775	&0.9673&0.9998&1.0000  \\

 \hline
  &  \multicolumn{2}{c}{$exp(0.2)$}  &  \multicolumn{2}{c}{$gamma(3,2)$}  &  \multicolumn{2}{c}{$Weibull(3,2)$}  &  \multicolumn{2}{c}{$lognormal(0.5,0.5)$}  &  \multicolumn{2}{c}{$Makeham(1,0.5)$} &  \multicolumn{2}{c}{$LFR(0.2)$} \\

\cline{1-13}
10&  0.0123&0.0593&0.2558&0.5651&0.8643	&0.9795
&0.4827&	0.7827&	0.1213	&0.3148	&0.3327&0.6508    \\
20& 0.0116&	0.0558&	0.6887&	0.9006&0.9998&	1.0000
&0.8969&0.9783&0.3400&0.6143	&0.8173	&0.9585    \\
30&0.0114&0.0557&0.9157&0.9846&	1.0000 &1.0000&0.9883&0.9981		&0.5770	&0.8169	&0.9693&0.9966   \\
40 & 0.0112	&0.0528	&0.9832	&0.9982	&1.0000 &1.0000
&0.9990	&0.9998		&0.7531&0.9178&0.9953&0.9997 \\
50& 0.0097&0.0496&0.9976&0.9998	&1.0000&	1.0000
&1.0000	&1.0000	&0.8724	&0.9692	&0.9995&1.0000  \\

 \hline
  &  \multicolumn{2}{c}{$exp(5)$}  &  \multicolumn{2}{c}{$gamma(2,1.5)$}  &  \multicolumn{2}{c}{$Weibull(2,3)$}  &  \multicolumn{2}{c}{$lognormal(3,0.8)$}  &  \multicolumn{2}{c}{$Makeham(2,0.5)$} &  \multicolumn{2}{c}{$LFR(3)$} \\

\cline{1-13}
10&  0.0133	&0.0612	&0.0814	&0.2569	&0.3394&0.6542&	0.1026&	0.1723
&0.2520&0.5229	&0.3354	&0.6598    \\
20& 0.0118&0.0602&0.2345&0.4975	&0.8157&0.9603&	0.1319&	0.2569
&0.6543	&0.8677&0.8105&0.9589   \\
30&0.0114&0.0561&0.4153&0.6929&	0.9714&	0.9966&	0.1566	&0.3383
&0.8768	&0.9700	&0.9709	&0.9972   \\
40 & 0.0107&0.0552&	0.5818	&0.8229&		0.9964&	0.9999	&	0.2125&	0.4124 	&0.9673&0.9947&0.9973&0.9998  \\
50& 0.0105&0.0503&0.7233&0.9074	&0.9998	&0.9999		&0.2632	&0.4729 & 0.9930&0.9996&0.9998&1.0000  \\
 \hline
\end{tabular}\end{adjustbox}
\end{table}

\begin{table}[h]
\caption{Power of the proposed test for various distributions-II }
\begin{adjustbox}{width=\textwidth}
\begin{tabular}{cccccccccccccccccc}\hline
\multirow{2}{*}{} & \multicolumn{12}{c}{$s=2$}\\ 
\cline{1-13}
$n/\alpha$ &  $0.01$  &  $0.05$ &  $0.01$  &  $0.05$ &  $0.01$  &  $0.05$ &  $0.01$  &  $0.05$ &  $0.01$  &  $0.05$ &  $0.01$  &  $0.05$  \\ \hline
 &  \multicolumn{2}{c}{$exp(1)$}  &  \multicolumn{2}{c}{$gamma(2,1)$}  &  \multicolumn{2}{c}{$Weibull(2,1)$}  &  \multicolumn{2}{c}{$lognormal(1,0.5)$}  &  \multicolumn{2}{c}{$Makeham(1,1)$} &  \multicolumn{2}{c}{$LFR(1)$}
 \\ \cline{1-13}

10& 0.0161&0.0622&0.2253&0.3816	&0.5486	&0.7181	&0.8424	&0.9416	&0.1034	&0.2068	&0.5271	&0.7142     \\
20& 0.0138&0.0616&0.4275&0.6289	&0.8419	&0.9485	&0.9958&	1.0000&0.1672&0.2995&0.8735	&0.9556   \\
30&  0.0132	&0.0552&0.6114&	0.7952&	0.9697&	0.9905&1.0000&	1.0000&0.2107&0.3736&0.9683&0.9921    \\
40& 0.0123&	0.0523&0.7596&0.8943&0.9968	&0.9996	&1.0000&	1.0000&0.2562&0.4295&0.9951&0.9994    \\
50&  0.0110	&0.0511&0.8580&	0.9476&	0.9999&	1.0000&1.0000&	1.0000&0.3023&0.4808&1.0000	&1.0000     \\
 \hline
 &  \multicolumn{2}{c}{$exp(2)$}  &  \multicolumn{2}{c}{$gamma(2,2)$}  &  \multicolumn{2}{c}{$Weibull(2,2)$}  &  \multicolumn{2}{c}{$lognormal(2,0.5)$}  &  \multicolumn{2}{c}{$Makeham(2,2)$} &  \multicolumn{2}{c}{$LFR(2)$} \\

\cline{1-13}
10& 0.0143& 0.0631& 0.2544& 0.4336& 0.5234& 0.7283& 0.8356& 	0.9509& 0.1296& 0.2156& 0.5209& 0.6867     \\
20& 0.0126	&0.0616	&0.4389&0.6301&	0.8541&0.9456&	0.9942	&1.0000&0.1477&	0.2648&	0.8523&	0.9332     \\
30& 0.0121&	0.0603&	0.5775&	0.7753&	0.9763&	0.9934&	1.0000	&1.0000&0.1894&	0.3390&	0.9661&	0.9931   \\
40&  0.0119&0.0514&0.7685&	0.8961&	0.9920&	0.9995&	1.0000&	1.0000&0.2525&0.4178&0.9942&0.9995    \\
50& 0.0105&0.0492&0.8841&0.9653	&1.0000	&1.0000	&1.0000	&1.0000	&0.3421	&0.5252	&0.9998&1.0000   \\
 \hline
 &  \multicolumn{2}{c}{$exp(0.5)$}  &  \multicolumn{2}{c}{$gamma(3,1)$}  &  \multicolumn{2}{c}{$Weibull(3,1)$}  &  \multicolumn{2}{c}{$lognormal(3,0.5)$}  &  \multicolumn{2}{c}{$Makeham(2,1)$} &  \multicolumn{2}{c}{$LFR(0.5)$} \\

\cline{1-13}
10& 0.0162&	0.0574&	0.5345&	0.7213&0.9357&0.9713&0.8492&	0.9510&	0.1618&0.3228&0.5209&	0.7196   \\
20&  0.0124&0.0565&	0.8747&	0.9572&0.9995&	1.0000&0.9991	&1.0000	&0.3069	&0.4757	&0.8593&0.9395    \\
30& 0.0123&0.0524&0.9692&0.9917	&1.0000	&1.0000	&1.0000&	1.0000&0.4467&0.6261&0.9773	&0.9962   \\
40& 0.0114&	0.0489&0.9982&0.9997&1.0000	&1.0000	&1.0000&	1.0000	&0.5628&0.7464	&0.9971&1.0000  \\
50& 0.0099&	0.0495&0.9961&1.0000&1.0000	&1.0000	&1.0000&	1.0000&	0.6674&0.8025&0.9998&1.0000    \\

 \hline
  &  \multicolumn{2}{c}{$exp(0.2)$}  &  \multicolumn{2}{c}{$gamma(3,2)$}  &  \multicolumn{2}{c}{$Weibull(3,2)$}  &  \multicolumn{2}{c}{$lognormal(0.5,0.5)$}  &  \multicolumn{2}{c}{$Makeham(1,0.5)$} &  \multicolumn{2}{c}{$LFR(0.2)$} \\

\cline{1-13}
10&  0.0157&0.0564	&0.5326&0.7134&	0.9384&0.9781&	0.8523&	0.9534&	0.2128&	0.3321&	0.5623&0.7361    \\
20& 0.0144	&0.0546&0.8509&	0.9556	&0.9992&1.0000&0.9957	&1.0000	&0.3197&0.4863&0.8413&0.9475     \\
30& 0.0136&0.0538&0.9764&0.9957	&1.0000	&1.0000	&1.0000	&1.0000	&0.4485&0.6434&	0.9709&	0.9912   \\
40 &  0.0122&0.0516&0.9975&	0.9999&1.0000&1.0000&1.0000	&1.0000	&0.5643&0.7361&	0.9978&	0.9995  \\
50& 0.0117&0.0512&	0.9999&	1.0000&	1.0000&1.0000&	1.0000	&1.0000	&0.6495	&0.8062&1.0000&	1.0000    \\

 \hline
  &  \multicolumn{2}{c}{$exp(5)$}  &  \multicolumn{2}{c}{$gamma(2,1.5)$}  &  \multicolumn{2}{c}{$Weibull(2,3)$}  &  \multicolumn{2}{c}{$lognormal(3,0.8)$}  &  \multicolumn{2}{c}{$Makeham(2,0.5)$} &  \multicolumn{2}{c}{$LFR(3)$} \\

\cline{1-13}
10& 0.0156&	0.0579&	0.2235&	0.3976&	0.5348&0.7182&0.1946	&0.3712	&0.3158&0.4616&0.5127&0.7083       \\
20& 0.0138&	0.0565&0.4156&0.6254&0.8655	&0.9557	&0.3258&	0.5154&0.5368&0.6856&0.8427&0.9337   \\
30& 0.0128&0.0526&0.6146&0.8028	&0.9695&0.9958&	0.5331	&0.7390&0.7113&	0.8475&0.9778&0.9948   \\
40 & 0.0112	&0.0516	&0.7571&0.8913	&0.9981&0.9992&	0.6317&	0.8309&	0.8424&0.9235	&0.9938&0.9993    \\
50&  0.0109& 0.0509	& 0.8542& 0.9489& 	0.9997& 1.0000	& 	0.7954& 0.9156	& 0.9413& 0.9786& 1.0000& 1.0000   \\
 \hline
\end{tabular}\end{adjustbox}
\end{table}

\begin{table}[h]
\caption{Power of the proposed test for various distributions-III }
\begin{adjustbox}{width=\textwidth}
\begin{tabular}{cccccccccccccccccc}\hline
\multirow{2}{*}{} & \multicolumn{12}{c}{$s=3$}\\ 
\cline{1-13}
$n/\alpha$ &  $0.01$  &  $0.05$ &  $0.01$  &  $0.05$ &  $0.01$  &  $0.05$ &  $0.01$  &  $0.05$ &  $0.01$  &  $0.05$ &  $0.01$  &  $0.05$  \\ \hline
 &  \multicolumn{2}{c}{$exp(1)$}  &  \multicolumn{2}{c}{$gamma(2,1)$}  &  \multicolumn{2}{c}{$Weibull(2,1)$}  &  \multicolumn{2}{c}{$lognormal(1,0.5)$}  &  \multicolumn{2}{c}{$Makeham(1,1)$} &  \multicolumn{2}{c}{$LFR(1)$}
 \\ \cline{1-13}

10& 0.0151& 0.0564& 0.3275& 0.4856& 0.5824& 0.7276&0.9234& 0.9843& 0.1617& 0.2474& 0.6035& 0.7265 \\
20& 0.0127&	0.0548&0.4716&0.641&0.8477&0.9182&1.0000&1.0000	&	0.1648&0.2536&0.8486&0.9257   \\
30& 0.0125&0.0526&0.6425&0.8165	&0.9557&0.9836&1.0000&1.0000	&0.1521	&0.2746	&0.9616&0.9877     \\
40& 0.0117&	0.0517&	0.7672&0.8876&0.9881&0.9942	&1.0000&	1.0000&	0.1923&	0.3145&	0.9886&0.9983    \\
50& 0.0112&0.0515&0.8638&0.9423	&0.9986&1.0000	&1.0000&	1.0000&0.2484&0.4056&0.9966&1.0000  \\
 \hline
 &  \multicolumn{2}{c}{$exp(2)$}  &  \multicolumn{2}{c}{$gamma(2,2)$}  &  \multicolumn{2}{c}{$Weibull(2,2)$}  &  \multicolumn{2}{c}{$lognormal(2,0.5)$}  &  \multicolumn{2}{c}{$Makeham(2,2)$} &  \multicolumn{2}{c}{$LFR(2)$} \\

\cline{1-13}
10&0.0151&0.0599&0.2955	&0.4286	&0.5985&0.7395&0.9017&0.9709	&0.1523	&0.2634	&0.6296&0.7393       \\
20& 0.0142&	0.0561&	0.5095&	0.6693	&0.8481&0.9216&	0.9998&	0.9996&0.1356&0.2480&0.8677&0.9434     \\
30& 0.0123&0.0547&0.6535&0.8124&0.9526&0.9845&	1.0000&	1.0000	&0.1675&0.2921&0.9673&0.9954    \\
40& 0.0117& 0.0523& 0.7735& 0.8863& 0.9854& 0.9965& 1.0000& 1.0000& 0.1736& 0.3115& 0.9871& 0.9974   \\
50& 0.0112&0.0516&0.8714&0.9584	&0.9976&0.9999&1.0000	&1.0000	&0.2147	&0.3746	&0.9982&0.9997  \\
 \hline
 &  \multicolumn{2}{c}{$exp(0.5)$}  &  \multicolumn{2}{c}{$gamma(3,1)$}  &  \multicolumn{2}{c}{$Weibull(3,1)$}  &  \multicolumn{2}{c}{$lognormal(3,0.5)$}  &  \multicolumn{2}{c}{$Makeham(2,1)$} &  \multicolumn{2}{c}{$LFR(0.5)$} \\

\cline{1-13}
10&  0.0152	&0.0593	&0.6438	&0.7785	&0.9287	&0.9645	&0.8923&	0.9662&0.2039&0.3141&0.6042&0.7473   \\
20&  0.0136	&0.0574	&0.8925	&0.9623	&0.9981	&1.0000	&1.0000	&1.0000	&0.3026&0.4264&0.8274&0.9128     \\
30&  0.0121	&0.0562	&0.9823	&0.9932	&1.0000	&1.0000	&1.0000	&1.0000	&0.3454	&0.4768	&0.9591&0.9968   \\
40 & 0.0118&0.0548&	0.9991&1.0000&	1.0000&	1.0000	&1.0000&	1.0000&0.4268&0.5730&0.9952	&1.0000   \\
50&  0.0112	&0.0567&0.9974&	1.0000	&1.0000&1.0000	&1.0000	&1.0000	&0.5080&0.6729&	0.9993	&1.0000    \\

 \hline
  &  \multicolumn{2}{c}{$exp(0.2)$}  &  \multicolumn{2}{c}{$gamma(3,2)$}  &  \multicolumn{2}{c}{$Weibull(3,2)$}  &  \multicolumn{2}{c}{$lognormal(0.5,0.5)$}  &  \multicolumn{2}{c}{$Makeham(1,0.5)$} &  \multicolumn{2}{c}{$LFR(0.2)$} \\

\cline{1-13}
10&  0.0144	&0.0612&0.6566	&0.7862	&0.9178&0.9616&	0.9157	&0.9734&0.2478	&0.3322	&0.5978&0.7324       \\
20&  0.0129	&0.0572	&0.8953	&0.9675&1.0000&1.0000&0.9997	&1.0000	&0.2623&0.3941	&0.8413&0.9155   \\
30& 0.0124& 0.0563	& 0.9751& 0.9942& 1.0000& 1.0000& 1.0000& 	1.0000& 0.3385& 0.4814& 0.9591& 0.9842   \\
40 & 0.0118&0.0544&0.9968&	0.9999	&1.0000	&1.0000	&1.0000	&1.0000	&0.4391&0.5865&0.9854&0.9972    \\
50& 0.0112	&0.0527&0.9998&1.0000&1.0000&1.0000	&1.0000	&1.0000	&0.4794&0.6347&1.0000&1.0000    \\

 \hline
  &  \multicolumn{2}{c}{$exp(5)$}  &  \multicolumn{2}{c}{$gamma(2,1.5)$}  &  \multicolumn{2}{c}{$Weibull(2,3)$}  &  \multicolumn{2}{c}{$lognormal(3,0.8)$}  &  \multicolumn{2}{c}{$Makeham(2,0.5)$} &  \multicolumn{2}{c}{$LFR(3)$} \\

\cline{1-13}
10&0.0153& 0.0558& 	0.2871& 0.4569& 0.5438	& 0.6799& 0.2946& 	0.4581& 0.3436& 0.4681	& 0.5835& 0.7087         \\
20& 0.0139&	0.0547&	0.4555	&0.6488&0.8455&0.9232&0.4995	&0.6956&0.4724&0.6153&0.8309&0.9183    \\
30& 0.0129	&0.0536	&0.6038&0.7872&0.9557&0.9856&0.7044&	0.8714&0.5987&0.7365&0.9625	&0.9854   \\
40 & 0.0116&0.0527&0.7759&0.9091&0.9862	&0.9988	&0.8336	&0.9418	&0.7041&0.8152&	0.9825&0.9978     \\
50& 0.0115&	0.0511&	0.8460&	0.9345&	0.9989&	1.0000		&0.9116&0.9777&	0.7986&	0.9024&0.9996&1.0000     \\
 \hline
\end{tabular}\end{adjustbox}
\end{table}
\end{document}